\documentclass[12pt]{article}

\usepackage{amssymb}
\usepackage{amsmath}
\usepackage{graphicx}
\usepackage{color}

\newtheorem{theorem}{Theorem}[section]
\newtheorem{corollary}[theorem]{Corollary}
\newtheorem{lemma}[theorem]{Lemma}
\newtheorem{proposition}[theorem]{Proposition}

\newtheorem{remark} [theorem]{Remark}

\newtheorem{example}[theorem]{Example}

\usepackage[shortlabels]{enumitem}
\newcommand{\qed}{\hfill{\rule{1ex}{1ex}}}
\DeclareMathOperator{\ind}{ind}

\DeclareMathOperator{\diag}{diag}

\DeclareMathOperator{\im}{im}
\DeclareMathOperator{\coker}{coker}
\newcommand{\vp}{\varphi}

\newcommand{\nn}{\nonumber}

\newcommand{\bA}{{\bf A}}
\newcommand{\bB}{{\bf B}}

\newcommand{\bD}{{\bf D}}

\newcommand{\bF}{{\bf F}}

\newcommand{\bK}{{\bf K}}

\newcommand{\cF}{\mathcal{F}}

\newcommand{\cJ}{\mathcal{J}}

\newcommand{\cL}{\mathcal{L}}

\newcommand{\cP}{\mathcal{P}}
\newcommand{\cQ}{\mathcal{Q}}
\newcommand{\cR}{\mathcal{R}}

\newcommand{\fB}{{\mathfrak B}}

\newcommand{\fN}{{\mathfrak N}}

\newcommand{\fS}{{\mathfrak S}}

\newcommand{\sC}{{\mathbb C}}

\newcommand{\sN}{{\mathbb N}}

\newcommand{\sR}{{\mathbb R}}

\newcommand{\sW}{{\mathbb W}}

\newcommand{\sZ}{{\mathbb Z}}

\begin{document}

\vspace*{10mm}

\begin{center}
{\Large\textbf{Wiener-Hopf plus Hankel operators: Invertibility Problems}}\footnote{This work was  supported by
the Special Project on High-Performance
Computing of the National Key R\&D Program of China (Grant No.~2016YFB0200604), the National Natural Science Foundation of China (Grant No.~11731006) and the Science Challenge Project of China (Grant No.~TZ2018001)..}
\end{center}

\vspace{5mm}

\begin{center}

\textbf{Victor D. Didenko and Bernd Silbermann}

\vspace{2mm}

Southern University of Science and Technology, Shenzhen, China; diviol@gmail.com

Technische Universit{\"a}t Chemnitz, Fakult{\"a}t f\"ur Mathematik,
09107 Chemnitz, Germany; silbermn@mathematik.tu-chemnitz.de

 \end{center}

  \vspace{10mm}



\vspace{18mm} \setcounter{page}{1} \thispagestyle{empty}

\begin{abstract}
The invertibility of Wiener-Hopf plus Hankel operators $W(a)+H(b)$ acting on the spaces
$L^p(\mathbb{R}^+)$, $1 < p<\infty$ is studied.  If $a$ and $b$ belong to a subalgebra of $L^\infty(\mathbb{R})$ and satisfy the condition
 \begin{equation*}
a(t) a(-t)=b(t) b(-t),\quad t\in\mathbb{R},
 \end{equation*}
we establish necessary and also sufficient conditions for the operators
$W(a)+H(b)$ to be one-sided invertible, invertible or generalized invertible. Besides,
efficient representations for the corresponding inverses are given.
\end{abstract}

%\tableofcontents % for long articles

\section{Introduction\label{s1}}

Let $\sR^-$ and $\sR^+$ be, respectively, the subsets of all
negative and all positive real numbers and $\chi_E$ refer to
the characteristic function of the subset $E\subset \sR$  -- i.e.
 \begin{equation*}
\chi_E(t):=\left \{
\begin{array}{ll}
  1 & \text{ if } t\in E,  \\
  0 & \text{ if } t\in \sR\setminus E  \\
\end{array}
 \right.
\end{equation*}
 In what follows, we often identify the spaces $L^p(\sR^+)$ and
$L^p(\sR^-)$, $1\leq p \leq\infty$ with the subspaces $\chi_{\sR^+} L^p(\sR)$ and $\chi_{\sR^-} L^p(\sR)$ of the space $L^p(\sR)$, which consist of the functions vanishing on $\sR^-$ and $\sR^+$, respectively.

Let $\mathcal{F}$ and $\mathcal{F}^{-1}$ be the direct and inverse Fourier transforms -- i.e.
\begin{equation*}
 \mathcal{F}\varphi(\xi):=\int\limits_{-\infty}^\infty e^{i\xi x} \varphi(x)\,dx,\quad
\mathcal{F}^{-1}\psi(x):=\frac1{2\pi}\int\limits_{-\infty}^\infty
e^{-i\xi x}\psi(\xi)\,d\xi,\;\;x\in\sR.
 \end{equation*}
Consider the set $\cL$ of functions $c\colon \sR\to \sC$ such that
$c=\cF k$ with $k\in L^1(\sR)$, and let $AP_W(\sR)\subset
L^\infty(\sR)$ be the set of functions $a\colon\sR\to \sC$ having the
representation
 \begin{equation}\label{Eq18}
 a(t)=\sum_{j\in\sZ} a_j e^{i\delta_j t}, \quad t\in\sR,
 \end{equation}
with absolutely convergent series \eqref{Eq18}. It is assumed that
that $\delta_j\in\sR$ for all $j\in \sZ$ and $\delta_j \neq \delta_k$ if $j\neq k$. By $G$ we denote the
smallest closed subalgebra of $L^\infty(\sR)$, which contains both
$AP_W(\sR)\subset
L^\infty(\sR)$ and $\cL$. One can show that any function $g\in G$ can be
represented in the form
  \begin{equation}\label{Eq19}
  g=a+c, \quad a\in AP_W(\sR), \; c\in \cL.
 \end{equation}
We also consider the subalgebra $G^+$  ($G^-$) of the algebra $G$,
which consists of all functions \eqref{Eq19} such that all numbers
$\delta_j$ are non-negative (non-positive) and the functions $c=\cF k$ such that $k(t)=0$  for all $t\leq 0$ ($t\geq 0$).  The functions from $G^+$ and
$G^-$ admit holomorphic extensions respectively to the upper and lower
half-planes and the set $G^+\cap G^-$ contains constant functions only.

Any function $a\in G$ generates an operator $W^0(a)\colon L^p(\sR)\to L^p(\sR)$ and
operators $W(a),H(a)\colon L^p(\sR^+)\to L^p(\sR^+)$ defined by
 \begin{equation*}%\label{TandH}
 \begin{aligned}
 W^0(a)&:=\mathcal{F}^{-1}a\mathcal{F}\varphi,  \\
  W(a)&:=PW^0(a) , \\
  H(a)&:= PW^0(a)QJ,
\end{aligned}
\end{equation*}
where $P \colon f\to \chi_{\sR^+} f$ and $Q:=I-P$ are the projections on
the subspaces $L^p(\sR^+)$ and $L^p(\sR^-)$, correspondingly, and
 $J\colon L^p(\sR)\to L^p(\sR)$ is the reflection operator defined by  $J\vp :=
\tilde{\vp}$. Here and in what follows, $\tilde{\vp}(t):=\vp(-t)$ for any $\vp\in L^p(\sR)$, $p\in [1,\infty]$. The operator
$W(a)$ is called the convolution on the semi-axis $\sR^+$ or the
Wiener-Hopf operator, whereas $H(a)$ is referred to as the Hankel
operator.  It is well-known \cite{GF:1974} that for $a\in G$ all three operators are bounded on the corresponding space $L^p$ for any  $p\in [1,\infty)$.

The operators $W^0$ and $W(a)$ can be also represented as
  \begin{align*}
 W^0(a)\vp(t)&  =\sum_{j=-\infty}^\infty a_j \vp(t-\delta_j)
 +\int_{-\infty}^\infty k(t-s) \vp(s)\,ds, \quad t\in\sR ,\\
  W(a)\vp(t)  &=\sum_{j=-\infty}^\infty a_j B_{\delta_j} \vp(t)
 +\int_{0}^\infty k(t-s) \vp(s)\,ds, \quad t\in\sR^+ ,
 \end{align*}
where
 \begin{align*}
B_{\delta_j} \vp(t)&=\vp(t-\delta_j) \quad\text{if}\quad \delta_j\leq 0,\\
B_{\delta_j} \vp(t)&=\left \{
 \begin{array}{l}
 0, \quad 0\leq t \leq  \delta_j\\
 \vp(t-\delta_j), \quad t> \delta_j
 \end{array}
 \right.
  \quad \text{if}\quad \delta_j >0.
 \end{align*}
Moreover, for $a=\mathcal{F}k$ the operator $H(a)$ acts as
\begin{equation*}
 H(a)\vp(t)=\int_0^\infty k(t+s)\vp(s)\,ds
\end{equation*}
and for  $a=e^{\delta t}$ as
 \begin{align*}
 H(a)\vp(t) &=
     \left \{
   \begin{array}{ll}
    \vp(\delta-t),& \quad 0\leq t\leq \delta\\
    0,& \quad t>\delta
     \end{array}
    \right .,
    \quad \text{if}\quad \delta >0,\\
   H(a)\vp(t)& =  0, \quad t\in \sR^+, \quad \text{if}\quad \delta \leq0.
 \end{align*}
Let us now recall a few useful identities involving the  operators mentioned. It is easily seen that if  $a,b\in G$, then
 \begin{equation*}
W^0(a b)=W^0(a) W^0(b).
\end{equation*}
Wiener-Hopf operators $W(a)$ generally do not possess this  property, but according to \cite[pp. 484, 485]{BS:2006} we still have
 \begin{equation}\label{cst4}
 \begin{aligned}
    W(ab)&=W(a)W(b)+H(a) H(\tilde{b}),\\
    H(ab)&= W(a)H(b)+H(a)W(\tilde{b}).
\end{aligned}
 \end{equation}
Moreover, if  $b\in G$, $c\in G^+$ and $c\in G^-$, then
 \begin{equation*}%\label{prod}
W(abc)=W(a) W(b) W(c).
 \end{equation*}

The operators $W(a)$ are well studied. For various classes of generating functions $a$, the conditions of Fredholmness or semi-Fredholmness of such operators can be efficiently written  \cite{BS:2006,BKS:2002,CD:1969,Du:1973,Du:1977,Du:1979,GF:1974}. Moreover, Fredholm and semi-Fredholm Wiener-Hopf operators are one-sided invertible,   the corresponding one-sided inverses are known and there is an efficient description of the kernels and cokernels of $W(a)$, $a\in G$.

Now we consider the Wiener-Hopf plus Hankel operators $\sW(a,b)$ acting on the space $L^p(\sR^+)$ and defined by
 \begin{equation}\label{WHH}
 \sW(a,b):=W(a)+H(b), \quad a,b\in L^\infty(\sR).
 \end{equation}
The study of such operators is much more involved. Nevertheless, Fredholm properties of \eqref{WHH} can be established either directly or by passing to a Wiener-Hopf operator with a matrix symbol. Thus Roch \emph{et al.}~\cite{RSS:2011} studied the Fredholmness of Wiener-Hopf plus Hankel operators with piecewise continuous generating functions, acting on $L^p$-spaces, $p\in [1,\infty)$. Another approach, called equivalence after extension, has been applied to operators with generating functions from a variety of classes. Nevertheless,  in spite of a vast amount of publications, it is mainly applied to the operators of a special form, namely, to the operators $\sW(a,a)=W(a)+H(a)$ acting on the $L^2$-space. It turns out that the Fredholmness, one-sided invertibility or invertibility of such operators  are equivalent to the corresponding properties of the Wiener-Hopf operator $W (a \tilde{a}^{-1})$,  so that they can be studied. However, even if an operator $\sW(a,a)$ is invertible, the corresponding inverse is not given (see, e.g. \cite[Corollary 2.2]{CN:2009} for typical results obtained by the method mentioned). If $a\neq b$, then hardly verified assumptions concerning the factorization of auxiliary matrix-functions is used. To study  Wiener-Hopf plus Hankel operators of the form $I+H(b)$, another method has been employed in \cite{KS:2000, KS:2001}, where the essential spectrum and the index of such operators are determined.

On the other hand, recently the Wiener-Hopf plus Hankel operators \eqref{WHH} have been studied under the assumption that the generating functions $a$ and $b$ satisfy the condition
  \begin{equation}\label{eqn1}
  a\tilde{a}=b\tilde{b}.
 \end{equation}
Thus if $a,b\in G$, then the Coburn-Simonenko Theorem for some classes of operators $\sW(a,b)$ is established \cite{DS:2014b}, and an efficient description of the space $\ker \sW(a,b)$ is obtained \cite{DS:2017a}. The aim of this work is to find conditions for one-sided invertibility, invertibility and generalized invertibility of the operators $\sW(a,b)$ and to provide efficient representations for the corresponding inverses when generating functions $a$ and $b$ satisfy the matching con\-di\-tion~\eqref{eqn1}. Similar problems for Toeplitz plus Hankel operators have been recently discussed in \cite{BE:2013, BE:2017, DS:2014a, DS:2014, DS:2016a, DS:2017}. However, the situation with Wiener-Hopf plus Hankel operators has some special features. Thus the operators here can also be semi-Fredholm -- i.e. in general, they may have infinite kernels and co-kernels. This creates additional difficulties. Therefore, in some cases, the results obtained are not as complete as for Fredholm Toeplitz plus Hankel operators.

This paper is organized as follows. Section~\ref{s2} contains known results on properties of Wiener-Hopf operators, Wiener-Hopf factorization of
functions $g\in G$ such that $g(t)g(-t)=1$, $t\in \sR$ and demonstrates their role in the description of the kernels of Wiener-Hopf plus Hankel operators $W(a)+H(b)$ under the condition \eqref{eqn1}.
In Section \ref{s3}, we establish necessary conditions for one-sided invertibility of the operators $\sW(a,b)$. Section~\ref{sec4} provides sufficient conditions for one-sided invertibility and presents efficient representations for the corresponding inverses. In Section \ref{s5}, we construct generalized inverses for Wiener-Hopf plus Hankel operators. The invertibility conditions presented in Section \ref{sec6} are supported by simple examples.

\section{Auxiliary results\label{s2}}

Let us recall the properties of Wiener-Hopf and Wiener-Hopf plus Hankel operators with generating functions from the algebra $G$. Thus it was shown in \cite{GF:1974} that for invertible functions $g$ the operators $W(g)$ are one-sided invertible. More precisely, if $a\in AP_W$, $c\in\cL$ and $g=a+c$ is invertible in $G$, then the element $a$ is also invertible in $G$.  Therefore, the numbers
 \begin{align}\label{ind}
    \nu(g):=\lim_{l\to\infty}\frac{1}{2l} [\arg g(t)]_{-l}^l,
\quad n(g):=\frac{1}{2\pi} [\arg
(1+g^{-1}(t)c(t)]_{t=-\infty}^\infty,
\end{align}
are  correctly defined. Moreover, the function $g$ admits the
factorization of the form
 \begin{equation}\label{Eq20}
 g(t)=g_-(t) e^{i\nu t}\left (\frac{t-i}{t+i} \right )^n g_+(t), \quad -\infty
 <t <\infty,
 \end{equation}
where $g_+^{\pm1}\in G^+$, $g_-^{\pm1}\in G^-$, $\nu=\nu(g)$ and
$n=n(g)$.

Let $-\infty<\nu<\infty$ be a real number.  On the space $L^p(\sR^+)$ we consider an operator $U_\nu$ defined by
  $$
(U_\nu \vp)(t):= \left \{
\begin{array}{ll}
  \vp(t-\nu) & \;\text{ if }\; \max(\nu,0) <t, \\[1ex]
  0  & \;\text{ if }\; 0\leq t\leq\max(\nu,0).
\end{array}
 \right .
  $$
It is easily seen that for any $\nu\geq0$, the operator $U_{\nu}$ is left invertible and $U_{-\nu}$ is one of its left-inverses. Moreover, $U_\nu=W(e^{it\nu})$   and $I-U_\nu U_{-\nu}$ is the projection operator
 $$
(I-U_\nu U_{-\nu} \vp)(t):= \left \{
\begin{array}{ll}
  \vp(t) & \;\text{ if }\; 0< t < \nu, \\[1ex]
  0  & \;\text{ if }\; \nu< t<\infty.
\end{array}
 \right .
  $$
We also consider operators $V$ and $V^{(-1)}$ defined by
 \begin{align*}%\label{eq}
(V\vp)(t):=\vp(t)-2\int_{0}^t e^{s-t}\vp(s)\,ds,  \,
(V^{(-1)}\vp)(t):=\vp(t)-2\int_t^{\infty} e^{t-s}\vp(s)\,ds.
 \end{align*}
 Set  $ V^{(m)}=V^m$ if $m\geq 0$ and
$V^{(-m)}=(V^{(-1)})^{-m}$ if $m<0$. It is known that if $m\in\sN$,
then $V^{(-m)} V^{(m)}=I$ \cite[Chapter 7]{GF:1974}, so that for $m>0$ the
operator $P_m:=I- V^{(m)}V^{(-m)}$ is a projection.

The factorization \eqref{Eq20} has been used in the construction of one-sided inverses for the Wiener-Hopf operators $W(g)$.

  \begin{theorem}[{\cite{GF:1974}}]\label{thm5} If $g=a+c\in G$, $a\in AP_W, c\in
\cL$, then the operator $W(g)$ is one-sided invertible in
$L^p(\sR^+)$, $1\leq p<\infty$ if and only if $g$ is invertible
in $G$. Moreover, if $g\in G$ is invertible in $G$ and $\nu:=\nu(g)$, $n:=n(g)$, then

\begin{enumerate}[(i)]

    \item  If $\nu>0$ and $n\geq0$, then the operator $W(g)$ is left invertible
and
  \begin{equation}\label{Eq21}
  W_l^{-1}(g)=W(g_+^{-1})V^{(-n)} U_{-\nu}  W(g_-^{-1})
 \end{equation}
is one of its left-inverses.
   \item  If $\nu>0$ and $n<0$, then the operator $W(a)$ is left invertible
and
  \begin{equation}\label{Eq22}
  W_l^{-1}(g)=W(g_+^{-1})(I-U_{-\nu} P_{-n} U_{\nu} )^{-1} U_{-\nu} V^{-n}
  W(g_-^{-1})
 \end{equation}
is one of its left-inverses.
  \item  If $\nu<0$ and $n\leq0$, then the operator $W(a)$ is right invertible and
  \begin{equation}\label{Eq23}
  W_r^{-1}(g)=W(g_+^{-1}) V^{-n} U_{-\nu}  W(g_-^{-1})
 \end{equation}
is one of its right-inverses.
   \item  If $\nu<0$ and $n>0$, then the operator $W(a)$ is right invertible and one of its right-inverses is
  \begin{equation}\label{Eq24}
  W_r^{-1}(g)=W(g_+^{-1}) V^{(-n)}U_{-\nu}(I-U_{\nu} P_{n} U_{-\nu} )^{-1}
  W(g_-^{-1}),
 \end{equation}
where
 \begin{equation}\label{Eq25}
 (I-U_{-\nu} P_{-n} U_{\nu} )^{-1}=\sum_{j=0}^\infty (U_{-\nu} P_{-n}
U_{\nu})^j,
 \end{equation}
and the series in the right-hand side of \eqref{Eq25} is uniformly
convergent.
\item  If $\nu=0$ and $n\leq0$ $(n\geq0)$, then the operator $W(g)$ is right (left) invertible and one of the corresponding inverses has the form
  \begin{equation}\label{Eq24a}
  W_{r/l}^{-1}(g)=W(g_+^{-1}) V^{(-n)} W(g_-^{-1}),
 \end{equation}
\end{enumerate}
  \end{theorem}
Let us point out that there is also an efficient description of the kernels of the operators $W(g)$, but the structure of $\ker W(g)$ depends on the indices $\nu(g)$ and $n(g)$ and it will be reminded later on.

As far as the Wiener-Hopf plus Hankel operators $\sW(a,b):=W(a)+H(b)$ are concerned, here we always assume that the generating functions $a,b\in G$  and satisfy the matching condition \eqref{eqn1}. In this case, the duo $(a,b)$ is referred to as the matching pair. Moreover, in what follows, we will only consider the matching pairs $(a,b)$ with the elements $a$ invertible in $G$. Notice that if $\sW(a,b)$ is semi-Fredholm, then $a$ is invertible in $G$ and the matching condition yields the invertibility of $b$ in $G$.

Let us introduce another pair $(c,d)$ with the elements $c$ and $d$
defined by
 \begin{equation*}%\label{Eq0}
  c:=ab^{-1}=  \tilde{a}^{-1}\tilde{b}, \quad d:=a\tilde{b}^{-1}=
  \tilde{a}^{-1}b.
 \end{equation*}
This duo is called the subordinated pair for $(a,b)$.The functions $c$ and $d$ possess a number of remarkable properties -- e.g. $c\tilde{c}=1=d\tilde{d}$. Following \cite{DS:2014b}, any function $g\in L_\infty(\sR)$ satisfying the condition $g \tilde{g}=1$
is called
 matching function. In passing note,  that if $(c,d)$ is
the subordinated pair for a matching pair $(a,b)$, then
$(\bar{d},\bar{c})$ is the subordinated pair for the
matching pair $(\bar{a}, \overline{\tilde{b}})$,  which
defines the adjoint operator
  \begin{equation}\label{cst10}
\sW^*(a,b)=W(\bar{a})+H(\overline{\tilde{b}}),
 \end{equation}
for the operator $\sW(a,b)$.

The next proposition comprises results from \cite{DS:2014b,DS:2017a}. For the reader's convenience, they are reformulated in a form suitable for subsequent presentation.

 \begin{proposition}\label{p1}
Assume that $g\in G$ is a matching function -- i.e. $g\tilde{g}=1$. Then
 \begin{enumerate}[(i)]
\item  Under the condition $g_-(0)=1$, the factors $g_+$ and $g_-$ in the factorization~\eqref{Eq20} are uniquely defined -- viz. the factorization takes the form
\begin{equation}\label{cst20}
g(t)=\left(\boldsymbol\sigma(g)\,\tilde{g}_+^{-1}(t)\right)e^{i\nu t}\Big(\frac{t-i}{t+i}\Big)^n g_+(t),
\end{equation}
where $\nu=\nu(g), n=n(g)$, $\boldsymbol\sigma(g)=(-1)^n g(0)$, $\tilde{g}_+^{\pm1}(t)\in G^-$ and $g_-(t)=\boldsymbol\sigma(g)\,\tilde{g}_+^{-1}(t)$.

\item If $\nu<0$ or if $\nu=0$ and $n<0$, then $W(g)$ is right-invertible and the operators $\mathbf{P}_g^{\pm}$,
   \begin{equation*}
\mathbf{P}_g^{\pm}:=(1/2)(I\pm JQВW^0(g)P) \colon \ker W(g)\to \ker W(g),
   \end{equation*}
considered on the kernel of the operator $W(g)$
are complementary projections.

\item If $(c,d)$ is the subordinated pair for a matching pair $(a,b)\in G \times G$ such that the operator $W(c)$ is right-invertible and $W_r^{-1}(c)$ is any right-inverse of $W(c)$, then
    \begin{align*}
\vp^+=\vp^+(a,b)&:=\frac{1}{2} (
W_r^{-1}(c)W(\tilde{a}^{-1})- JQW^0(c)P
W_r^{-1}(c)W(\tilde{a}^{-1}))\nn\\
&\quad + \frac{1}{2} JQW^0(\tilde{a}^{-1}),
  \end{align*}
is an injective operator from  $\ker W(d)$ into $\ker (W(a)+H(b))$.

 \item\label{iv} If $(c,d)$ is the subordinated pair for the matching pair $(a,b)$, then

\begin{enumerate}[(a)]

  \item If the operator $W(c)\colon  L^p(\sR^+)\to L^p(\sR^+)$, $1<p<\infty$ is right-invertible, then
  \begin{equation}\label{ker}
\ker(W(a)+H(b))=\vp^+(\im \mathbf{P}^+_d) \dotplus\im \mathbf{P}^-_c.
  \end{equation}
  \item If the operator $W(d)\colon L^p(\sR^+)\to L^p(\sR^+)$, $1<p<\infty$ is left-invertible, then
\begin{equation}\label{coker}
\coker(W(a)+H(b))=\vp^+(\im \mathbf{P}^+_{\bar{c}}) \dotplus\im
\mathbf{P}^-_{\bar{d}},
\end{equation}
where the operator $\vp^+$ in \eqref{coker} is defined by the matching pair $(\bar{a}, \bar{\widetilde{b}})$.
 \end{enumerate}

 \item\label{v} Let $\Lambda_j$ be the normalized Laguerre polynomials and the functions $\psi_j$, $j\in \sZ_+$, be defined by
 \begin{align}\label{lag}
    \psi_j(t)&:= \left \{
\begin{array}{ll}
\sqrt{2} e^{-t} \Lambda_j(2t),& \text{ if } t>0,\\
 0, & \text{ if } t<0,\\
 \end{array}
     \right .,& \quad j=0,1,\cdots\,.\phantom{--}
     \end{align}
Then for $\nu=0$ and $n<0$,
the following systems $\fB_{\pm}(g)$ of functions $W(g_+^{-1})
\psi_{j}$ form bases in the spaces $\im \mathbf{P}^{\pm}_g$:
 \begin{enumerate}[(a)]
\item If $n=-2m$, $m\in\sN$, then
 \begin{equation*}
\fB_{\pm}(g)=\{W(g_+^{-1}) \left ( \psi_{m-k-1}\mp
\boldsymbol\sigma(g)\psi_{m+k}\right ): k=0,1,\cdots, m-1\},
 \end{equation*}
 and
  \begin{equation}\label{even}
  \dim\im \mathbf{P}^{\pm}_g=m.
  \end{equation}
 \item If $n=-2m-1$, $m\in\sZ_+$, then
   \begin{equation*}
   \fB_{\pm}(g)=\{W(g_+^{-1})\left ( \psi_{m+k}\mp
\boldsymbol\sigma(g)\psi_{m-k}\right ): k=0,1,\cdots, m\}\setminus \{0\},
   \end{equation*}
 and
  \begin{equation}\label{odd}
  \dim\im \mathbf{P}^{\pm}_g=m+ \frac{1\mp\boldsymbol\sigma(g)}{2}.
  \end{equation}
   \end{enumerate}

     \end{enumerate}
   \end{proposition}
\begin{remark}
If $\nu<0$, the corresponding spaces $\im \mathbf{P}^{\pm}_g$ are also described in~\cite{DS:2017a}. However, these representations are not used in what follows so that they are not included to the above proposition.
\end{remark}

\section{Necessary conditions for one-sided\\ inver\-ti\-bi\-lity\label{s3}}
From now on we always assume without mentioning it specifically  that the generating functions $a$ and $b$  constitute a matching pair. Moreover, let us also recall that if an operator  $W(a)+H(b)$, $a,b\in G$ acting in the space $L^p(\sR^+)$, $p\in (1,\infty)$ is Fredholm or semi-Fredholm, then the generating function $a$ is invertible in $G$. Therefore, the elements $c$ and $d$ of the subordinated pair $(c,d)$ are also invertible in $G$ and the Wiener-Hopf operators $W(c)$ and $W(d)$ are Fredholm or semi-Fredholm. Let $\nu_1:=\nu(c)$, $n_1:=n(c)$, $\nu_2:=\nu(d)$, and  $n_2:=n(d)$ be the corresponding indices \eqref{ind} of the functions $c$ and $d$. We start with necessary conditions for one-sided invertibility of the operators $W(a)+H(b)$ in the case where at least one of the indices $\nu_1$, $\nu_2$ is not equal to zero. The situation $\nu_1=\nu_2=0$ will be considered later on.

 \begin{theorem}\label{thm3.1}
 Let $a,b\in G$ and the operator $W(a)+H(b)$ be one-sided invertible in $L^p(\sR^+)$. Then:
  \begin{enumerate}[(i)]
    \item  Either $\nu_1 \nu_2\geq 0$ or $\nu_1>0$ and $\nu_2<0$.
    \item If $\nu_1=0$ and $\nu_2>0$, then  $n_1 >-1$ or $n_1=-1$ and $\boldsymbol\sigma(c)=-1$.
    \item  If $\nu_1<0$ and $\nu_2=0$, then $n_2<1$ or $n_2=1$ and $\boldsymbol\sigma(d)=-1$.
  \end{enumerate}
  \end{theorem}

\textit{Proof}
 \emph{(i)} Assume that $\nu_1 \nu_2 < 0$. If $\nu_1<0$ and $\nu_2>0$, then the operator $W(c)$ is right invertible whereas $W(d)$ is left invertible. Moreover, the kernel of the operator $W(c)$ and cokernel of $W(d)$ are infinite-dimensional \cite{GF:1974} and so are the spaces $\im \mathbf{P}^{-}_c$ and $\im \mathbf{P}^{-}_{\overline{d}}$ \cite[Theorems 2.4 and 2.5]{DS:2017a}. Taking into account Proposition~\ref{p1}\ref{iv}, we obtain that $\ker (W(a)+H(b))\neq \{0\}$ and $\coker (W(a)+H(b))\neq \{0\}$, hence the operator $W(a)+H(b)$ is not one-sided invertible.

\emph{(ii)} Let $\nu_2>0$. By Proposition~\ref{p1}\ref{iv}, the operator $W(a)+H(b)$ has a non-zero cokernel. If, in addition, $n_1<-1$ or $n_1=1$ and $\boldsymbol\sigma(c)=1$, then \eqref{even} and  \eqref{odd} show that in both cases, $\im \mathbf{P}^-_c \neq\{0\}$. Therefore, according to \eqref{ker}, the operator $W(a)+H(b)$ also has a non-trivial kernel and is not one-sided invertible.

The assertion (iii) can be proved analogously.
\qed

 Let us briefly discuss the case where $\nu_1>0$ and $\nu_2<0$. As was mentioned in \cite{DS:2017a}, in this situation it is not clear whether the corresponding Wiener-Hopf operator is even normally solvable. Nevertheless, the kernel and cokernel of $W(a)+H(b)$ can still be described. This gives a possibility to establish necessary conditions of one-sided invertibility. However, they are not as transparent as before and, in addition to the relations between the indices $\nu_1, \nu_2, n_1,n_2$, the corresponding conditions can include information about the factors in the Wiener-Hopf factorizations of the subordinated functions $c$ and $d$.  We consider one of possible cases.

  \begin{theorem}\label{thm3.2}
 Let $\nu_1>0$, $\nu_2<0$, $n_1=n_2=0$ and let $\fN_\nu^p$, $\nu>0$ denote the set of functions $f\in L^p(\sR^+)$ such that $f(t)=0$ for $t\in (0,\nu)$.

 \begin{enumerate}[(i)]
   \item If the operator $W(a)+H(b)\colon L^p(\sR^+)\to L^p(\sR^+)$, $1<p<\infty$ is invertible from the left, then
     \begin{equation}\label{leftinv}
  \vp^+(\mathbf{P}^+_d)\cap \fN_{\nu_1/2}^p= \{0\},
     \end{equation}
  where $\vp^+=\vp^+(ae^{-i\nu_1 t/2}, be^{i\nu_1 t/2})$.

   \item If the operator $W(a)+H(b)\colon L^p(\sR^+)\to L^p(\sR^+)$, $1<p<\infty$ is invertible from the right, then
     \begin{equation}\label{rightinv}
  \vp^+(\mathbf{P}^+_{\overline{c}})\cap \fN_{-\nu_2/2}^p= \{0\},
     \end{equation}
where $\vp^+=\vp^+(\overline{a}e^{i\nu_2 t/2},\overline{\tilde{b}}e^{-i\nu_2 t/2})$.
 \end{enumerate}
   \end{theorem}
   \textit{Proof} Let $\nu_1>0$, $\nu_2<0$, $n_1=n_2=0$   and  $W(a)+H(b)$ be a left-invertible operator. It can be represented in the form
   \begin{equation}\label{eqn3.1}
W(a)+ H(b)= \left ( W \left ( ae^{-i\nu_1 t/2} \right )+ H\left
(
 be^{i\nu_1 t/2} \right )  \right ) W \left ( e^{i\nu_1 t/2} \right
 ).
   \end{equation}
Direct computations show that $( ae^{-i\nu_1 t/2}, be^{i\nu_1 t/2})$  is a matching pair with the subordinated pair $(c_1,d_1)=(c e^{-i\nu_1 t}, d)$. Since $\nu(c_1)=0$, $n(c_1)=n_1=0$, the kernel of the operator $W(c_1)$ is trivial. Consequently, $\ker P^-_{c_1}=\{0\}$ and the relation \eqref{ker} yields
 \begin{equation*}
\ker\left ( W \left ( ae^{-i\nu_1 t/2} \right )+ H\left
(
 be^{i\nu_1 t/2} \right )  \right )=\vp^+(\im \mathbf{P}^+_d)
 \end{equation*}
with the operator $\vp^+=\vp^+(ae^{-i\nu_1 t/2}, be^{i\nu_1 t/2})$.
Therefore, taking into account \eqref{eqn3.1}, we  obtain
 \begin{equation*}
\ker (W(a)+H(b)) =\{\eta=W(e^{-i\nu_1t/2})u: u\in \vp^+(\mathbf{P}^+_d)\cap \im W(e^{i\nu_1t/2}) \}.
 \end{equation*}
If the operator $W(a)+H(b)$ is left invertible, its kernel consists of the zero element only. However, since $\im W(e^{i\nu_1t/2})=\fN_{\nu_1/2}^p$ and
  \begin{equation*}
 \ker W(e^{-i\nu_1t/2})\cap (\vp^+(\mathbf{P}^+_d)\cap \fN_{\nu_1/2}^p)=\{0\},
\end{equation*}
the assumption
   \begin{equation*}
  \vp^+(\mathbf{P}^+_d)\cap \fN_{\nu_1/2}^p \neq \{0\}
\end{equation*}
yields the non-triviality of the kernel of $W(a)+H(b)$, so that \eqref{leftinv} holds.

The second assertion in Theorem \ref{thm3.2} comes from the first one by passing to the adjoint operator (see \eqref{cst10}).
   \qed

   \begin{remark}\label{rem0}
Theorem \ref{thm3.2} raises an interesting question: Are there exist invertible operators $W(a)+H(b)$, such that
\begin{equation*}
\dim \coker W(c)=\dim\ker W(d)=\infty?
\end{equation*}
Note that in the case $\nu(c)=\nu(d)=0$, for any prescribed natural number $N$ one can probably find invertible operators $W(a)+H(b)$ for which
  \begin{equation}\label{large}
 \ind |W(c)| > N, \quad  \ind |W(d)| >N.
\end{equation}
Note that the set of Toeplitz plus Hankel operators possesses the property~\eqref{large} -- cf. \cite{DS:2019a}, but for Wiener-Hopf plus Hankel operators, this problem requires a separate study.
   \end{remark}

  \begin{remark}\label{rem1}
Although the description of the spaces $\im \mathbf{P}^+_d$ and $\im \mathbf{P}^+_{\overline{c}}$ is available \cite{DS:2017a}, the verification of the conditions \eqref{leftinv}-\eqref{rightinv} is not trivial. It depends on the properties of Wiener-Hopf operators constituting the operator $\vp^+$ and may require a lot of effort.
  \end{remark}

\begin{remark}\label{rem2}
If $\nu_1>0$, $\nu_2<0$ but $n_1\neq 0$ or/and $n_2\neq 0$, the necessary conditions of one-sided invertibility have the same form \eqref{leftinv} and \eqref{rightinv} but the representation \eqref{eqn3.1},  spaces $\fN_\nu^p$ and  operators $\vp^+$ should be redefined accordingly.
  \end{remark}

We now consider the situation when both indices $\nu_1$ and $\nu_2$ vanish.
Let us start with an auxiliary result.

 \begin{lemma}\label{lem2}
 If $(a,b)\in G\times G$ is a matching pair with the subordinated pair $(c,d)$, then for the factorization signatures of the functions $c$ and $d$ the equation
 \begin{equation}\label{sign}
\boldsymbol\sigma(c)=\boldsymbol\sigma(d)
 \end{equation}
holds and  the indices $n_1$ and $n_2$ are simultaneously odd or even.
 \end{lemma}
\textit{Proof}
Let $n(a)$ and $n(b)$ be the corresponding indices \eqref{ind} for the functions $a$ and $b$, respectively. Then
\begin{equation}\label{ccc}
n_1=n(c)=n(a b^{-1})=n(a)-n(b), \quad n_2=n(d)=n(a \tilde{b}^{-1})=n(a)+n(b).
\end{equation}
Therefore,
\begin{align*}%\label{eqn}
  \boldsymbol\sigma(c) & = (-1)^{n(a)-n(b)}c(0)= (-1)^{n(a)-n(b)}a(0)b^{-1}(0),  \\
  \boldsymbol\sigma(d) & = (-1)^{n(a)+n(b)}d(0)= (-1)^{n(a)+n(b)}a(0)\tilde{b}^{-1}(0),
\end{align*}
and since $b(0)=\tilde{b}(0)$ and the numbers $n(a)-n(b)$ and $n(a)+n(b)$ are simultaneously odd or even, the equation \eqref{sign} follows.

 Moreover, using the relations \eqref{ccc} again, we obtain
 \begin{equation*}
n_1+n_2=2 n(a),
\end{equation*}
so that $n_1$ has the same evenness as $n_2$.
\qed

 We start with the left invertibility of the operators $\sW(a,b)$.

\begin{theorem}\label{thm3.3}
 If $a,b\in G$, $\nu_1=\nu_2=0$, $n_2\geq n_1$ and the operator $W(a)+H(b)$ is invertible from the left, then the index $n_1$ satisfies the inequality
     \begin{equation*}
     n_1\geq -1
     \end{equation*}
     and if $n_1=-1$, then $\boldsymbol\sigma(c)=-1$ and $n_2>n_1$.

  \end{theorem}
  \textit{Proof}
If $n_1<-1$, then the operator $W(c)$ is right invertible. By Proposition~\ref{p1}\ref{v}, the image of the projection $\mathbf{P}^-_c$ contains non-zero elements, and by \eqref{ker} so is $\ker (W(a)+H(b))$. This contradicts the left invertibility of the operator $W(a)+H(b)$, hence $n_1\geq -1$.

Assume now that $n_1=-1$ and $\boldsymbol\sigma(c)=1$. Using \eqref{ker} and Proposition~\ref{p1}\ref{v} again, we note that $\im \mathbf{P}_c^-\neq \{0\}$, so that the operator $W(a)+H(b)$ has a non-trivial kernel and, therefore, it is not left-invertible. Hence $\boldsymbol\sigma(c)=-1$. Assuming, in addition,  that $n_2=-1$ and  $\ker (W(a)+H(b))=\{0\}$, we obtain
\begin{equation*}
\boldsymbol\sigma(c)=-1, \quad \boldsymbol\sigma(d)=1,
\end{equation*}
which is not possible by  Lemma \ref{lem2}. Hence, $n_2>n_1$.
 \qed

 \begin{theorem}\label{thm3.4}
If $a,b\in G$, $\nu_1=\nu_2=0$, $n_1>n_2$ and the operator $W(a)+H(b)$ is invertible from the left, then the inequality
 \begin{equation*}%\label{nec}
 n_1\geq 1
 \end{equation*}
holds. Moreover, the index $n_2$ is either non-negative or $n_2<0$ and $n_1\geq -n_2$.
                  \end{theorem}
   \textit{Proof}
If $n_1\leq0$, then $n_2\leq -2$ -- cf. Lemma~\ref{lem2}, and  $W(a)+H(b)$ has a non-trivial kernel, which contradicts the left-invertibility of this operator.  On the other hand, if $1\leq n_1$ and  $0\leq n_2$, then $W(a)+H(b)$ is clearly left-invertible, so we proceed with the case $n_2<0$.  By Lemma~\ref{lem2}, both numbers $n_1$ and $n_2$ are either even or odd.  In both cases the proof of the fact that the indices $n_1$ and $n_2$ satisfy the inequality $n_1\geq n_2$ is similar, but each situation should be examined separately. Here we only analyse the case where  $n_1$ and $n_2$ are odd numbers. Considering $\ind W(c):=\mathbf{k}_1=-n_1$ we chose  $k_1\in \sZ$ such that
\begin{equation*}%\label{indices}
 2k_1 + \mathbf{k}_1=1.
\end{equation*}
Then, according to \cite[Theorem 3.2]{DS:2017a}, we have
\begin{equation}\label{eqnKer}
 \begin{aligned}
  & \ker(W(a)+H(b)) = \left \{ W\left ( \left ( \frac{t-i}{t+i} \right)^{-k_1}\right )u :\right .  \\
    & \left . u\in \left \{  \frac{1+\boldsymbol\sigma(c)}{2}W(c_+^{-1})\{\sC\psi_0\} \dotplus \vp^+(\im \mathbf{P}^+_d) \right \} \cap \im W\left (\! \left ( \!\frac{t-i}{t+i} \right)^{k_1}\!\right )\!\right \},
 \end{aligned}
\end{equation}
where the operator $\vp^+=\vp^+(a_1,b_1)$ is defined by the matching pair
\begin{equation*}
(a_1,b_1)=\left (a(t)\left(\frac{t-i}{t+i} \right)^{-k_1}, b(t)\left(\frac{t-i}{t+i} \right)^{k_1}\right )
\end{equation*}
and $c_+$ is the plus factor in the Wiener-Hopf factorization \eqref{cst20} of the function $c$. The function $\psi_0$ is defined in \eqref{lag} and using another representation of the Laguerre polynomials -- cf.~\cite[Eq.~(2.5)]{DS:2017a}, one can show that
 \begin{equation*}
 \im W\left ( \left ( \frac{t-i}{t+i} \right)^{k_1}\right ) = \mathrm{clos}\, \mathrm{span}_{L^p(\sR^+)} \left \{ \psi_{k_1}, \psi_{k_1+1}, \cdots \right \}.
 \end{equation*}
Thus if a function
 $$
 u \in \im W\left ( \left ( \frac{t-i}{t+i} \right)^{k_1}\right )
  $$
is expanded in a Fourier series of the Laguerre polynomials $\psi_j, j=0,1,\cdots$, its first $k_1$ Fourier-Laguerre coefficients are equal to zero. If we now assume that the dimension of the subspace
\begin{equation*}
\fS(c,d):= \left \{  \frac{1+\boldsymbol\sigma(c)}{2}W(c_+^{-1})\{\sC\psi_0\} \dotplus \vp^+(\im \mathbf{P}^+_d) \right \}
\end{equation*}
is greater than $k_1$, then there is a non-zero function $u_0\in \fS(c,d)$, the first $k_1$ Fourier-Laguerre coefficients of which vanish. Hence,  \eqref{eqnKer} shows that the kernel of $W(a)+H(b)$ contains a non-zero element. This contradicts the left invertibility of the operator $W(a)+H(b)$. Therefore,
 \begin{equation}\label{dimen}
k_1 \geq \dim \fS(c,d),
 \end{equation}
and taking into account the Eq.~\eqref{odd}, we rewrite the inequality \eqref{dimen} as
 \begin{equation}\label{aaa}
k_1 \geq \frac{1+\boldsymbol\sigma(c)}{2}+ k_2,
 \end{equation}
where
 \begin{equation*}
k_2= r +\frac{1-\boldsymbol\sigma(d)}{2}
 \end{equation*}
and $-n_2=2r+1$. Since $k_1=(1-\mathbf{k}_1)/2=(1+n_1)/2$,  the inequality \eqref{aaa} takes the form
 \begin{equation*}
\frac{1 + n_1}{2}\geq \frac{1 +\boldsymbol\sigma(c)}{2} + \frac{-n_2-1}{2} +\frac{1 -\boldsymbol\sigma(d)}{2}
 \end{equation*}
or
\begin{equation*}
n_1\geq -n_2 +\boldsymbol\sigma(c)-\boldsymbol\sigma(d).
\end{equation*}
Since $\boldsymbol\sigma(c)=\boldsymbol\sigma(d)$ by Lemma~\ref{lem2}, the proof is completed.
\qed

Thus Theorems \ref{thm3.3}, \ref{thm3.4} provide necessary conditions for the left invertibility of the  operators $\sW(a,b)=W(a)+H(b)$. Passing to right-invertible operators, one can recall a simple fact that the operator $\sW(a,b)$ is right-invertible if and only if the adjoint operator $\sW^*(a,b)$ is left invertible. However, relation \eqref{cst10} shows that
 \begin{equation*}
\sW^*(a,b)=\sW(\bar{a}, \overline{\tilde{b}}).
 \end{equation*}
We note that $(\bar{a}, \overline{\tilde{b}})$ is also a matching pair with the subordinated pair \linebreak $(c_1,d_1)=(\bar{d}, \bar{c})$, so that
\begin{equation*}%\label{bbb}
\begin{aligned}
 &\nu(c_1)=\nu(\bar{d})=-\nu_2, \quad\nu(d_1)= \nu(\bar{c})=-\nu_1,\\
   &n(c_1)= n(\bar{d})=-n_2, \quad n(d_1)= n(\bar{c})=-n_1,\\
   &\boldsymbol\sigma(c_1)=\boldsymbol\sigma(d), \quad \boldsymbol\sigma(d_1)=\boldsymbol\sigma(c).
\end{aligned}
\end{equation*}
Now Theorems \ref{thm3.3} and \ref{thm3.4} can be used to write the necessary conditions for the right invertibility of the operators $W(a)+H(b)$. Let us just formulate the corresponding results.

\begin{theorem}\label{thm3.6}
 Let $a,b\in G$, $\nu_1=\nu_2=0$, $n_1\leq n_2$ and the operator $W(a)+H(b)$ is invertible from the right. Then
 \begin{equation*}
 n_2\leq 1
 \end{equation*}
and if $n_2=1$, then $\boldsymbol\sigma(d)=-1$ and $n_1<n_2$.
            \end{theorem}

 \begin{theorem}\label{thm3.7}
Let $a,b\in G$, $\nu_1=\nu_2=0$, $n_1>n_2$ and the operator $W(a)+H(b)$ is invertible from the right. Then the inequality
    \begin{equation*}
  n_2\leq -1
 \end{equation*}
holds. Moreover,  the index $n_1$ is either non-positive or  $n_1\leq -n_2$.
                   \end{theorem}

\section{Sufficient conditions of one-sided\\ invertibility and one-sided inverses\label{sec4}}

Our next goal is to establish sufficient conditions for one-sided invertibility of the operators $W(a)+H(b)$. In fact, many necessary conditions above are also sufficient ones.
 \begin{theorem}\label{3.10} Let $a,b\in G$ and indices $\nu_1, \nu_2, n_1$ and $n_2$ satisfy any of the following conditions:
  \begin{enumerate}[(i)]
    \item $\nu_1<0$ and $\nu_2<0$.
    \item $\nu_1>0$, $\nu_2<0$, $n_1=n_2=0$, operator $W(a)+H(b)$ is normally solvable and satisfies the condition \eqref{rightinv}.
        \item $\nu_1<0$, $\nu_2=0$ and $n_2<1$ or $n_2=1$ and $\boldsymbol\sigma(d)=-1$.
    \item $\nu_1=0$, $n_1\leq0$ and $\nu_2<0$.
    \item $\nu_1=0$ and $\nu_2=0$
        \begin{enumerate}[(i)]
          \item $n_1\leq0$, $n_2<1$;
          \item $n_1\leq 0$,  $n_2=1$ and $\boldsymbol\sigma(d)=-1$;
                        \end{enumerate}
  \end{enumerate}
 Then the operator $W(a)+H(b)\colon L^p(\sR^+)\to L^p(\sR^+)$, $1<p<\infty$ is right invertible.
  \end{theorem}
  \textit{Proof}
By Proposition \ref{p1}, each condition in Theorem~\ref{3.10} ensures that
   \begin{equation*}
 \coker (W(a)+H(b))=\{0\}.
 \end{equation*}
Hence the operator $W(a)+H(b)$ is invertible from the right.
 \qed

Sufficient conditions for the left invertibility of the operators $W(a)+H(b)$ can be obtained from Theorem \ref{3.10} by passing to the  adjoint operators and we leave it to the reader.

In the remaining part of this section we deal with the construction of  right inverses for the operators $W(a)+H(b)$. Recall that one-sided inverses of Wiener-Hopf operators can be easily determined from the Wiener-Hopf factorizations of the corresponding generating functions -- cf.~Theorem~\ref{thm5}. However, finding the inverses for Wiener-Hopf plus Hankel operators is much more difficult problem and to the best of our knowledge, so far there was no efficient representation of the corresponding inverses even for the simplest pairs of generating functions. Now we want to establish formulas for the left and right inverses of the operators $W(a)+H(b)$ in the case of matching generating functions.

Let us assume that the operators $W(c)$ and $W(d)$ are invertible from the same side.  This condition is not necessary for the one-sided invertibility and note that the corresponding inverses can be also constructed even if the condition mentioned is not satisfied.
 \begin{theorem}\label{t5}
 Let $(a,b)\in  G \times G$ be a matching pair such that the
operators $W(c)$ and $W(d)$ are invertible from the right. Then the
operator $W(a)+H(b)$ is also right invertible and one of its right inverses has the form
  \begin{equation}\label{EqRI}
  B:= (I - H(\tilde{c})) \mathbf{A} + H(a^{-1})W_r^{-1}(d),
 \end{equation}
where $\mathbf{A}=W_r^{-1}(c) W(\tilde{a}^{-1})W_r^{-1}(d)$.
 \end{theorem}
  \textit{Proof}
The proof of this result uses equations \eqref{cst4}. Consider the product
$(W(a)+H(b))B$,
  \begin{equation}\label{EqN17}
  \begin{aligned}
&(W(a)+H(b))B  =\\
 &\quad (W(a)+H(b))(I - H(\tilde{c})) \mathbf{A} +
(W(a)+H(b))H(a^{-1})W_r^{-1}(d).
 \end{aligned}
 \end{equation}
It follows from \eqref{cst4} that
\begin{equation*}%\label{EqN21}
\begin{aligned}
H(b)H(\tilde{c})&=W(bc)-W(b)W(c)
     =W(a)-W(b)W(c), \\
      W(a)H(\tilde{c})&=H(a\tilde{c})-H(a)W(c)=H(b)-
      H(a)W(c).
 \end{aligned}
\end{equation*}
Therefore, the first product in the right-hand side of \eqref{EqN17} can be rewritten as
\begin{equation}
 \begin{aligned}\label{EqN3}
   &\quad(W(a)+H(b)) (I-H(\tilde{c})\mathbf{A}  =
   (W(b)W(c)+H(a)W(c))\mathbf{A} \\
   &=(W(b)\!+\!H(a))W(c)\mathbf{A}\!=\! (W(b)\!+\!H(a)\!)W(c)W_r^{-1}(c)
   W(\tilde{a}^{-1})W_r^{-1}(d) \\
   & =W(b)W(\tilde{a}^{-1})W_r^{-1}(d)+ H(a)\!
   W(\tilde{a}^{-1})W_r^{-1}(d).
   \end{aligned}
    \end{equation}
Analogously,
\begin{equation*}%\label{EqN4}
\begin{aligned}
    W(a)H(a^{-1})&=H(aa^{-1})-H(a)W(\tilde{a}^{-1})=-H(a)W(\tilde{a}^{-1}), \\
     H(b)H(a^{-1})&=W(b\tilde{a}^{-1})-W(b)W(\tilde{a}^{-1})
     =W(d)-W(b)W(\tilde{a}^{-1}),
\end{aligned}
 \end{equation*}
and the second product in the right-hand side of \eqref{EqN17} has the form
 \begin{equation}\label{EqN4}
 \begin{aligned}
&\quad (W(a)+H(b))H(a^{-1})W_r^{-1}(d)\\
&=-H(a)W(\tilde{a}^{-1})W_r^{-1}(d)
 +
W(d)W_r^{-1}(d)-W(b)W(\tilde{a}^{-1})W_r^{-1}(d)\\
 &=I-H(a)W(\tilde{a}^{-1})W_r^{-1}(d)-W(b)W(\tilde{a}^{-1})W_r^{-1}(d).
\end{aligned}
 \end{equation}
 Combining \eqref{EqN3} and \eqref{EqN4}, one obtains
 $$
(W(a)+H(b))B=I,
 $$
hence $B$ is a right inverse for the Wiener-Hopf plus
Hankel operator $W(a)+H(b)$.
  \qed

\begin{corollary}\label{cor2}
Let $a,b\in G$ and indices $\nu_1, \nu_2, n_1$ and $n_2$ satisfy one of the following conditions:
  \begin{enumerate}[(i)]
    \item $\nu_1<0$ and $\nu_2<0$.
    \item $\nu_1<0$, $\nu_2=0$ and $n_2\leq 0$.
    \item $\nu_1=0$, $n_1\leq0$ and $\nu_2<0$.
      \end{enumerate}
 Then $W(a)+H(b)$ is invertible from the right and one of its  right inverses can be  constructed by  formula \eqref{EqRI}.
 \end{corollary}

  \begin{example}
Let us consider the operator
  \begin{equation}\label{Eq26}
  \sW(\nu_1,\nu_2)=W(e^{i\nu_1t})+H(e^{i\nu_2t}),\quad t\in \sR,
 \end{equation}
where $\nu_1$ and $\nu_2$ are real numbers such that
    \begin{align}
 \nu_1-\nu_2\leq 0, \label{Eq27a}\\
    \nu_1+\nu_2\leq 0. \label{Eq27b}
 \end{align}
In passing note that the conditions \eqref{Eq27a}-\eqref{Eq27a} are equivalent to the inequality
\begin{equation*}
\nu_1\leq -|\nu_2|,
\end{equation*}
so that $\nu_1\leq0$. Consider now the generating functions
$a(t)=e^{i\nu_1t}$ and $b(t)=e^{i\nu_2t}$. They satisfy the matching
conditions \eqref{eqn1}, namely,
  $$
a(t)a(-t)=b(t)b(-t)=1.
  $$
The elements $c$ and $d$ of the subordinated pair for
 matching pair $(e^{i\nu_1t},e^{i\nu_2t})$ are
 $$
c(t)=e^{i(\nu_1-\nu_2)t}, \quad d(t)=e^{i(\nu_1+\nu_2)t}.
 $$
Taking into account the conditions  \eqref{Eq27a}-\eqref{Eq27a}, we observe that the corresponding Wiener-Hopf operators $W(c)$, $W(d)$ are  right invertible and have infinite dimensional kernels.
In order to construct a right inverse of the
operator \eqref{Eq26} one can use Theorem~\ref{t5}.
Let us recall simple properties of Wiener-Hopf
and Hankel operators with exponential generating function. Thus for
the generating function $a(t)=e^{i\nu t}$ one has:
 \begin{enumerate}[(i)]
    \item If $\nu \leq 0$, then the operator $W(e^{i\nu t})$ is
    right invertible and one of its right inverses is
      $$
W^{-1}_r (e^{i\nu t})=W(e^{-i\nu t}).
      $$
    \item If $\nu \geq 0$, then the operator $W(e^{i\nu t})$ is
    left invertible and one of its left inverses is
      $$
W^{-1}_l (e^{i\nu t})=W(e^{-i\nu t}).
      $$
    \item If $\nu<0$, then $H(e^{i\nu t})=0$.
 \end{enumerate}
Therefore,
  \begin{equation*}%\label{Eq28}
  W_r^{-1}(c)=W(e^{-i(\nu_1-\nu_2)t}), \quad
  W_r^{-1}(d)=W(e^{-i(\nu_1+\nu_2)t}).
 \end{equation*}
Thus the operator \eqref{Eq26} is subject to Theorem \ref{t5}. In order to write the cor\-responding right inverse of $W(a)+H(b)$, we first determine the operator~$\mathbf{A}$. Simple computations show that
 $$
\mathbf{A}=W(e^{-i(\nu_1-\nu_2)t}) W(e^{-i\nu_2t}).
 $$
Therefore the right inverse \eqref{EqRI} for the operator
\eqref{Eq26} has the form
 \begin{align*}
  (W(e^{i\nu_1t})+H(e^{i\nu_2t}))_r^{-1}=&
(I-H(e^{-i(\nu_1-\nu_2)t}))W(e^{-i(\nu_1-\nu_2)t})
W(e^{-i\nu_2t}) \\
 & + H(e^{-i\nu_1t})W(e^{-i(\nu_1+\nu_2)t})
 \nn.
\end{align*}
Moreover, using formulas  \eqref{cst4}, one obtains
 \begin{equation*}
H(e^{-i(\nu_1-\nu_2)t})W(e^{-i(\nu_1-\nu_2)t})=0, \quad H(e^{-i\nu_1t})W(e^{-i(\nu_1+\nu_2)t})= W(e^{-i\nu_1t}),
\end{equation*}
and the operator $(W(e^{i\nu_1t})+H(e^{i\nu_2t}))_r^{-1}$ can be
finally written as
 \begin{equation*}%\label{Eq29}
(W(e^{i\nu_1t})+H(e^{i\nu_2t}))_r^{-1}=
H(e^{-i\nu_1t})W(e^{-i(\nu_1+\nu_2)t}) +W(e^{-i\nu_1t}).
 \end{equation*}
   \end{example}

We now construct a left inverse for the operator $W(a)+H(b)$.

 \begin{theorem}\label{t4.4}
 Let $(a,b)\in  G \times G$ be a matching pair such that the
operators $W(c)$ and $W(d)$ are invertible from the left. Then the
operator $\sW(a,b)=W(a)+H(b)$ is also left-invertible and one of its left-inverses has the form
  \begin{equation}\label{EqLI}
  \sW_l(a,b)= \mathbf{C}(I - H(\tilde{d})) + W_l^{-1}(c)H(\tilde{a}^{-1}),
 \end{equation}
where $\mathbf{C}=W_l^{-1}(c) W(\tilde{a}^{-1})W_l^{-1}(d)$.
 \end{theorem}
  \textit{Proof}
Recalling that the adjoint operator $\sW^*(a,b)$ to the operator $W(a)+H(b)$ can be identified with the operator
\begin{equation*}%\label{adjoint}
\sW^*(a,b)=T(a_1)+H(b_1),\quad a_1=\bar{a}, \quad b_1=\overline{\tilde{b}},
\end{equation*}
we note that $(a_1,b_1)$ is a matching pair with the subordinated pair $(c_1,d_1)=(\bar{d}, \bar{c})$. Since $W(c_1)=W(\bar{d})$ and $W(d_1)=W(\bar{c})$ are invertible from the right, the operator $\sW^*(a,b)$ is also right-invertible by Theorem \ref{t5} and according to \eqref{EqRI}, one of its right inverses can be written as
 \begin{equation} \label{RInvW}
 \begin{aligned}
(\sW^*(a,b))_r^{-1}&=(I-H(\tilde{c}_1) \textbf{A}_1 +H(a_1^{-1})W_r^{-1}(d_1)\\
   &=(I-H(\overline{\tilde{d}}) \textbf{A}_1 +H(\bar{a}^{-1})W_r^{-1}(\bar{c}),
 \end{aligned}
 \end{equation}
where
\begin{equation}\label{A_1}
\textbf{A}_1=W_r^{-1}(c_1) W(\tilde{a}^{-1}_1) W_r^{-1}(d_1)=W_r^{-1}(\bar{d}) W(\overline{\tilde{a}}^{-1}) W_r^{-1}(\bar{c}).
\end{equation}
The left inverse to the operator $\sW(a,b)$ can be now obtained by computing the adjoint operator for the operator $(\sW^*(a,b))_r^{-1}$. Since for any right-invertible operator $A$ one has
 \begin{equation*}
\left ( A_r^{-1} \right )^* =\left ( A^* \right )_l^{-1},
 \end{equation*}
we can use the relations
 \begin{equation*}
W^*(g)=W(\bar{g}), \quad H^*(g)=H(\overline{\tilde{g}}), \quad g\in G,
 \end{equation*}
 to obtain the representation \eqref{EqLI} from \eqref{RInvW}-\eqref{A_1}.
     \qed

\section{Generalized invertibility of Wiener-Hopf\\ plus Hankel operators\label{s5}}

An operator $A$ is called generalized invertible if there exists an operator
$A_g^{-1}$ referred to as generalized inverse for $A$, such that
 $$
A A_g^{-1} A =A.
 $$
If $A_g^{-1}$ is a generalized inverse for the operator $A$ and
the equation
 \begin{equation}\label{EqG}
 Ax=y
 \end{equation}
is solvable, then the element $x_0=A_g^{-1} y$ is a solution of
equation \eqref{EqG}.

Our next task is to determine generalized inverses for Wiener-Hopf
plus Hankel operators $W(a)+H(b)$ if the generating functions $a$ and $b$ constitute a matching pair. For this we recall useful formulas connecting Wiener-Hopf plus Hankel operators and matrix Wiener-Hopf operators. Thus according to
\cite[Equation (2.4)]{DS:2014b}, the diagonal operator $\diag
(W(a)+H(b)+Q, W(a)-H(b)+Q)$ can be represented in the form
\begin{equation}\label{eqn3}
 \left(\!%
\begin{array}{cc}
  W(a)+H(b)+Q &  0\\
   0 &  W(a)-H(b)+Q\\
   \end{array}%
\!\right) =
 \cJ A_1 A_2(W(V(a,b))\!+\!\cQ)C \cJ^{-1}\,,
\end{equation}
where the operators $A_1,  A_2$, $\cJ$, $C=C(a,b)$ and $V=V(a,b)$
are defined by
 \begin{align}
 &  A_1:=\diag(I,I)-\diag(P,Q)W^0 \left(%
\begin{array}{cc}
  a & b \\
  \tilde{b}  & \tilde{a}  \\
   \end{array}%
\right) \diag(Q,P), \nn\\[1ex]   %\label{Eq12}\\
& A_2:=\diag(I,I)+\cP W^0 (V(a,b)) \cQ, \quad \cJ:=  \frac{1}{2}\left(%
\begin{array}{cc}
  I & J \\
 I & -J \\
\end{array}
\right) ,\nn\\
 & C(a,b) :=\left(%
\begin{array}{cc}
  I & 0 \\
 W^0(\tilde{b}) &  W^0(\tilde{a}) \\
   \end{array}%
\right),
 \quad
  V(a,b) :=\left(%
\begin{array}{cc}
   a-b \tilde{b} \tilde{a}^{-1} &  b \tilde{a}^{-1} \\
 - \tilde{b}\tilde{a}^{-1} & \tilde{a}^{-1} \\
   \end{array}%
\right), \nn
  \end{align}
and
 $$
\cP:=\diag(P,P), \quad \cQ:=\diag(Q,Q)\, .
  $$
Using  the notation
\begin{align}\label{Eq3.5}
B&:= W(V(a,b))+\cQ, \\[1ex]%\label{Eq14}\\
\cR&:= \diag(W(a)+H(b),W(a)-H(b)), \quad  R:= \cR+\cQ, \nn
\end{align}
we write the equation \eqref{eqn3} as
 \begin{equation}\label{Eq4}
 R= \cJ A_1 A_2 B C \cJ^{-1}\,.
 \end{equation}

Considering the operator $R$ and taking into account the equation \eqref{Eq4} and the invertibility of the operators $\cJ, C,
A_1$ and $A_2$, we write
 \begin{equation*}%\label{Eq4.5}
      R^{-1}_g =  \cJ C^{-1}  B^{-1}_g   A_2^{-1} A_1^{-1} \cJ^{-1}.
 \end{equation*}
Observe that $R^{-1}_g$ is diagonal since  so is the operator $R$.
Thus
 $$
R^{-1}_g =\diag(\bF^{-1}_{g},\bK^{-1}_{g}),
 $$
and it is clear that the diagonal elements $\bF^{-1}_{g}$ and
$\bK^{-1}_{g}$ have the form
 $$
\bF^{-1}_{g}=F^{-1}_{g}+ Q, \quad \bK^{-1}_{g}=K^{-1}_{g}+ Q,
 $$
where $F^{-1}_{g},K^{-1}_{g}\colon L^p(\sR^+)\to L^p(\sR^+)$ are
generalized inverses for the operators $W(a)+H(b)$ and $W(a)-H(b)$,
respectively.

In this section we construct a generalized inverse for  operator
$W(a)+H(b)$ provided that the operator $B$ is generalized invertible
and a generalized inverse of $B$ can be represented in a special
form. The following theorem has been proved in \cite{DS:2016a} in the case of Toeplitz plus Hankel operators. For Wiener-Hopf plus Hankel operators the proof literally repeats all constructions there and is omitted here.

 \begin{theorem}\label{t4}
Let $(a,b)$ be a matching pair with the subordinated pair $(c,d)$.
Assume that the operator $B$ of \eqref{Eq3.5} is generalized
invertible and has a generalized inverse $B^{-1}_g$ of the form
 \begin{equation}\label{Eq5}
  B^{-1}_g=
 \left(%
\begin{array}{cc}
 \bA  & \bB \\
   \bD & 0 \\
   \end{array}%
\right) +\cQ,
 \end{equation}
where  $\bA,\bB$ and $\bD$ are operators acting in the space
$L^p(\sR^+)$. Then $W(a)+H(b)$ is generalized invertible and the
operator $G$,
  \begin{equation}\label{Eq6}
\begin{aligned}
G: =-H(\tilde{c})(\bA(I-H(d))-\bB
H(\tilde{a}^{-1}))+H(a^{-1})\bD(I-H(d))+W(a^{-1}),
\end{aligned}
 \end{equation}
 is a generalized inverse for the operator $W(a)+H(b)$.
 \end{theorem}

\begin{lemma}\label{lem1}
Let $(a,b)\in G\times G$ be a matching pair such that one of the following conditions holds:
 \begin{enumerate}[(i)]
    \item The operators $W(c)$ and $W(d)$ are
    right invertible.
    \item The operators $W(c)$ and $W(d)$ are
    left invertible.
    \item $W(c)$  and $W(d)$ are, respectively,  left and right invertible
    operators.
 \end{enumerate}
Then the operator $B$ of \eqref{Eq3.5} is generalized invertible and
it has a generalized inverse of the form \eqref{Eq5}.
\end{lemma}
 \textit{Proof} For a matching pair $(a,b)$ the operator
$W(V(a,b))$ has the form
  $$
  W(V(a,b))=
  \left(%
\begin{array}{cc}
  0 & W(d) \\
 -W(c)   & W(\tilde{a}^{-1}) \\
   \end{array}%
\right).
  $$
Assume for definiteness that both operators $W(c)$ and $W(d)$ are
right invertible. Then the operator $W(V(a,b))$ is also right
invertible and it is easily seen that one of its right inverses is
given by the formula
 \begin{equation*}%\label{Eq14}
B_g^{-1}=
  \left(%
\begin{array}{cc}
 W_r^{-1}(c) W(\tilde{a}^{-1})W_r^{-1}(d)  &   -W_r^{-1}(c)\\[1ex]
 W_r^{-1}(d)   & 0 \\
   \end{array}%
\right) +\cQ,
 \end{equation*}
where $W_r^{-1}(c)$ and $W_r^{-1}(d)$ are right-inverses of the
operators $W(c)$ and $W(d)$, correspondingly. Thus in this case,
condition \eqref{Eq5} is satisfied with the operators
    \begin{equation}\label{Eq15}
\mathbf{A}= W_r^{-1}(c) W(\tilde{a}^{-1})W_r^{-1}(d), \quad
\mathbf{B}= -W_r^{-1}(c), \quad \mathbf{D}=W_r^{-1}(d).
 \end{equation}
The other cases are considered analogously. Thus if both
operators $W(c)$ and $W(d)$ are left invertible, then $B$ is left invertible with a left-inverse having the form \eqref{Eq5}, where
    \begin{equation}\label{Eq16}
\mathbf{A}= W_l^{-1}(c) W(\tilde{a}^{-1})W_l^{-1}(d), \quad
\mathbf{B}= -W_l^{-1}(c), \quad \mathbf{D}=W_l^{-1}(d),
 \end{equation}
and if $W(c)$ is left-invertible and $W(d)$ is right invertible,
then the corresponding operators $\mathbf{A}$, $\mathbf{B}$ and
$\mathbf{D}$ in \eqref{Eq5} are
     \begin{equation}\label{Eq17}
\mathbf{A}= W_r^{-1}(c) W(\tilde{a}^{-1})W_l^{-1}(d), \quad
\mathbf{B}= -W_r^{-1}(c), \quad \mathbf{D}=W_l^{-1}(d),
 \end{equation}
which completes the proof.
 \qed

Combining Theorem \ref{t4} and Lemma \ref{lem1} one obtains the
following result.

 \begin{theorem}\label{thm2.3}
Let operators $W(c)$ and $W(d)$ satisfy one of the assumptions
of Lemma \ref{lem1} and $\mathbf{A}, \mathbf{B}$ and
$\mathbf{D}$ be the operators defined by one of the
relations \eqref{Eq15}--\eqref{Eq17}. Then the operator $W(a)+H(b)$
is generalized invertible and \eqref{Eq6} is one of its generalized inverses.
 \end{theorem}

\begin{remark}
We note that in cases (i) and (ii), the operator $W(a)+H(b)$ is one-sided invertible and formulas for the corresponding inverses obtained in Section~\ref{sec4} are simpler that the representation \eqref{Eq6}.
 \end{remark}

\section{Invertibility of Wiener-Hopf plus Hankel operators\label{sec6}}

The results of the previous sections can now be used to establish various invertibility conditions for the operators $W(a)+H(b)$ and write the corresponding inverses. Let us formulate one of such results and provide a few examples.

 \begin{corollary}\label{cor1}
 Let $(a,b)$, $a,b\in  G $ be a matching pair such that the
operators $W(c)$ and $W(d)$ are invertible. Then the operator
$W(a)+H(b)$ is invertible and
  \begin{equation}\label{EqN5}
(W(a)+H(b))^{-1}=(I - H(\tilde{c})) W^{-1}(c)
W(\tilde{a}^{-1})W^{-1}(d) + H(a^{-1})W^{-1}(d).
 \end{equation}
  \end{corollary}

\textit{Proof}
If the operators $W(c)$ and $W(d)$ are invertible, then relations
(2.7) and (2.4) of \cite{DS:2014b} show that the operators
$W(a)+H(b)$ is invertible and the result follows from Theorem
\ref{t5}.
\qed

Let us point out that this is a very surprising result. There is a
vast literature devoted to the study of the Fredholmness and
one-sided invertibility of Wiener-Hopf plus Hankel operators in the
situation where generating functions satisfy the relation $b=a$
or $\tilde{b}=a$. Of course, such generating functions
constitute a matching pair. The other case studied is $a(t)=1$ for
all $t\in\sR$ and $b=b(t)$ is a specific matching function. However,
to the best of our knowledge, so far there are no efficient representations
for the inverse operators. On the other hand, for a wide class of generating
functions $g$ the inverse operators $W^{-1}(g)$ can be constructed.
Therefore, formula \eqref{EqN5} is an efficient
tool in constructing the inverse operators $(W(a)+H(b))^{-1}$ in the case where $a$ and $b$ constitute a matching generating pair.

\begin{example}\label{ex2}
Let us consider the operators $W(a)+H(b)$ in the case where $a=b$.
In this situation $c(t)=1$ and $d(t)=a(t) \tilde{a}^{-1}$.
Hence, $H(\tilde{c})=0$, $W(c)=I$ and if the operator $W(d)$ is
invertible, then the operator $W(a)+H(a)$ is also invertible and
 $$
(W(a)+H(a))^{-1}=(W(\tilde{a}^{-1}) + H(a^{-1}))W^{-1}(a
\tilde{a}^{-1}).
 $$
\end{example}

\begin{example}\label{ex3}
Let $b=\tilde{a}$. Then $c(t)=a(t)
\tilde{a}^{-1}(t)$ and $d(t)=1$. Hence, if the operator $W(c)$ is
invertible, then the operator $W(a)+H(a)$ is also invertible and
 $$
(W(a)+H(\tilde{a}))^{-1}=(I - H(\tilde{a}a^{-1})) W^{-1}(a
\tilde{a}^{-1}) W(\tilde{a}^{-1}) + H(a^{-1}).
 $$
\end{example}

\begin{example}\label{ex4}
Let $a(t)=1$  and $b(t)b(-t)=1$ for all $t\in
\sR$. In this situation, $c(t)=\tilde{b}(t)$, $d(t)=b(t)$  and
if the operator $T(b)$ is invertible, then
 $$
(I+H(b))^{-1}=(I - H(b)) W^{-1}(\tilde{b}) W^{-1}(b).
 $$
\end{example}

\section*{Conclusion}

For matching generating functions $a,b\in G$, the invertibility of the operators $W(a)+H(b)$ can be described in terms of indices $\nu$ and $n$ of the subordinated functions $c$ and $d$. Moreover, the corresponding inverses can be represented using only auxiliary  Wiener-Hopf and Hankel operators along with the corresponding inverses of scalar Wiener-Hopf operators.  This approach is efficient and can be realised as soon as the Wiener-Hopf factorization of the functions $c$ and $d$ is available --- cf.~\eqref{Eq21}-\eqref{Eq24a}.

 \end{document}